\documentclass[letterpaper, 10pt, conference]{ieeeconf}
\IEEEoverridecommandlockouts
\usepackage{amssymb}
\usepackage{amsmath}

\usepackage{amsthm}
\usepackage{xcolor}
\usepackage{cite}
\usepackage{hyperref}
\usepackage{tikz}
\usepackage{algorithm}
\usepackage{algpseudocode}

\usepackage{pifont} 
\newcommand{\cmark}{\ding{51}}
\newcommand{\xmark}{\ding{55}}

 \newcommand{\hquad}{\hspace{0.5em}} 
 
\newcommand{\argmin}[1]{\underset{{#1}}{\text{argmin}}}

\def\betastep{\beta_{\text{step}}}
\def\betamax{\beta_{\text{max}}}

\theoremstyle{definition}
\newtheorem{definition}{Definition}
\newtheorem{theorem}{Theorem}

\newtheorem{lemma}{Lemma}[theorem]

\def\HInf{\mathcal{H}_{\infty}}
\def\H2{\mathcal{H}_2}
\def\L1{\mathcal{L}_1}
\def\LInf{\mathcal{L}_{\infty}}
\def\D{\mathcal{D}}
\def\N{\mathcal{N}}
\def\P{\mathcal{P}}
\def\S{\mathcal{S}}
\def\RHInf{\mathcal{RH}_{\infty}}

\def\row{\text{row}}
\def\col{\text{col}}
\def\sub{\text{sub}}

\newcommand{\infinfnorm}[1]{\|{#1}\|_{\infty\rightarrow\infty}}
\newcommand{\oneonenorm}[1]{\|{#1}\|_{1\rightarrow 1}}
\newcommand{\oneinfnorm}[1]{\|{#1}\|_{1\rightarrow\infty}}
\newcommand{\infonenorm}[1]{\|{#1}\|_{\infty\rightarrow 1}}
\newcommand{\htwonorm}[1]{\|{#1}\|_{\H2}}
\newcommand{\perfnorm}[1]{\|{#1}\|_{\text{perf}}}
\newcommand{\stabnorm}[1]{\|{#1}\|_{\text{stab}}}
\newcommand{\frobnorm}[1]{\|{#1}\|_{F}}

\def\G{\mathbf{G}}
\def\K{\mathbf{K}}
\def\Lbf{\mathbf{L}}
\def\u{\mathbf{u}}
\def\v{\mathbf{v}}
\def\w{\mathbf{w}}
\def\x{\mathbf{x}}
\def\y{\mathbf{y}}
\def\z{\mathbf{z}}
\def\d{\hat{\mathbf{w}}} % Avoid confusion w/delta
\def\Phibf{\mathbf{\Phi}}
\def\Psibf{\mathbf{\Psi}}
\def\Lambdabf{\mathbf{\Lambda}}
\def\Phiset{\Phibf_{\{i,j\}}}
\def\Psiset{\Psibf_{\{i,j\}}}

\def\Phix{\mathbf{\Phi}_x}
\def\Phiu{\mathbf{\Phi}_u}
\def\Phixu{\begin{bmatrix} \Phix \\ \Phiu \end{bmatrix}}
\def\xu{\begin{bmatrix} \x \\ \u \end{bmatrix}}
\def\ZAB{\begin{bmatrix} zI-A & -B \end{bmatrix}}
\def\ZAC{\begin{bmatrix} zI-A \\ -C \end{bmatrix}}

\def\Phiw{\mathbf{\Phi}_w}
\def\Phiv{\mathbf{\Phi}_v}
\def\Phiwv{\begin{bmatrix} \Phiw & \Phiv \end{bmatrix}}
\def\wv{\begin{bmatrix} \w \\ \v \end{bmatrix}}

\def\Phixw{\mathbf{\Phi}_{xw}}
\def\Phixv{\mathbf{\Phi}_{xv}}
\def\Phiuw{\mathbf{\Phi}_{uw}}
\def\Phiuv{\mathbf{\Phi}_{uv}}
\def\Phiof{\begin{bmatrix} \Phixw & \Phixv \\ \Phiuw & \Phiuv \end{bmatrix}}
\def\IOHoriz{\begin{bmatrix} I & 0 \end{bmatrix}}
\def\IOVert{\begin{bmatrix} I \\ 0 \end{bmatrix}}

\title{\LARGE \bf
Distributed Robust Control for Systems with Structured Uncertainties}

\author{Jing Shuang (Lisa) Li and John C. Doyle
	\thanks{J. S. Li and J. C. Doyle are with Computing and Mathematical Sciences, California Institute of Technology. {\tt\small \{jsli, doyle\}@caltech.edu}. This research was in part supported by funding from NSERC PGSD3-557385-2021.}
}

\begin{document}

\maketitle

\begin{abstract}

We present D-$\Phi$ iteration: an algorithm for distributed, localized, and scalable robust control of systems with structured uncertainties. This algorithm combines the System Level Synthesis (SLS) parametrization for distributed control with stability criteria from $\L1$, $\LInf$, and $\nu$ robust control. We show in simulation that this algorithm achieves good nominal performance while greatly increasing the robust stability margin compared to the LQR controller. To the best of our knowledge, this is the first distributed and localized algorithm for structured robust control; furthermore, algorithm complexity depends only on the size of local neighborhoods and is independent of global system size. We additionally characterize the suitability of different robustness criteria for distributed and localized computation.

\end{abstract}
\section{Introduction}
Robust control theory \cite{Dahleh1993, Packard1993, Zhou1998} provides stability and performance guarantees in the face of uncertainties, which arise from modeling inaccuracies or unexpected operating conditions. Such uncertainties are common in engineering applications, and robust control theory plays a large role in the control and design of such systems. 
For large-scale systems, distributed $\HInf$ control has been studied for linear spatially invariant systems \cite{Bamieh2002} and linear symmetric systems \cite{Lidstrom2017}, and other special systems; more recently, distributed and localized synthesis for arbitrary linear systems using $\L1$ is proposed in \cite{Matni2020}. However, to the best of our knowledge, no distributed methods address \textit{structured} uncertainty, which often arises in large systems; for instance, in a power grid, we may have parametric uncertainty on electrical properties of transmission lines between buses, but no uncertainty between unconnected buses -- this imposes structure on system uncertainty. Incorporating knowledge of this structure into the synthesis procedure allows control engineers to provide less conservative robust stability margins, which has been historically important in the design of aerospace systems.

In this paper, we leverage the System Level Synthesis (SLS) parametrization \cite{Anderson2019} to provide distributed synthesis methods for structured robust control. We deal with diagonal time-varying uncertainty and use $\L1$, $\LInf$, and $\nu$ \cite{CompanionPaper_Nu} robustness criteria to formulate synthesis problems. We also present D-$\Phi$ iteration, an algorithm for distributed synthesis. To the best of our knowledge, this is the first distributed and localized algorithm for structured robust control.

We begin with a review of the SLS formulation for state feedback, full control, and output feedback in Section \ref{sec:sls}. We then describe the \textit{separability} of each formulation and how this impacts distributed computation in Section \ref{sec:separability}. In Section \ref{sec:robust_stability}, we describe how $\L1$, $\LInf$, and $\nu$ robustness criteria translated into separable objectives and constraints on the control problem, and present two versions of the scalable D-$\Phi$ iteration algorithm for robust control. The efficacy of the algorithms is demonstrated via simulations in Section \ref{sec:simulations}; D-$\Phi$ iteration provides good robust stability margins while mostly preserving nominal performance.

\section{System Level Synthesis} \label{sec:sls}

The SLS framework formulates control problems as optimizations over closed-loop responses $\Phibf$ \cite{Anderson2019}. We cast a control problem as an optimization of the form:
\begin{equation} \label{eq:sls_problem_generic}
    \min_{\Phibf} \quad f(\Phibf) \quad \textrm{s.t.} \quad \Phibf \in \S_a \cap \P
\end{equation}
where $f$ is any convex functional, and $\S_a$ and $\P$ are convex sets. $\S_a$ represents the \textit{achievability constraint}, which is enforced in all SLS problems; it ensures that we only search over closed-loop responses $\Phibf$ which can be achieved using a causal, internally stable controller. This is made mathematically explicit in Theorem \ref{thm:sf_main}. $\P$ represents additional, optional constraints; we typically use $\P$ to enforce sparsity on $\Phibf$, corresponding to local communication, delayed communication, or local disturbance rejection. Examples of how to formulate standard control problems (e.g. LQR, LQG) in the form of \eqref{eq:sls_problem_generic} are given in \cite{Anderson2019}; robust stability problems can also be cast in this form. For ease of numerical computation, we generally assume that $\Phibf$ is finite impulse response (FIR) with horizon $T$. We now provide definitions of $\Phibf$ and $\S_a$ for the state feedback, full control, and output feedback cases.

Consider the discrete-time linear time-invariant (LTI) system with state $\x$, control $\u$, and disturbance $\w$ in frequency domain (i.e. $z$-domain):
\begin{equation} \label{eq:plant_x}
    z\x = A\x + B\u + \w
\end{equation}
and apply linear causal state feedback controller $\u=\K\x$, where $\K$ is a transfer matrix and can be written in terms of spectral elements $K(p)$ as $\K = \sum_{p=1}^{\infty} K(p) z^{-p}$. Define $\Phix$ and $\Phiu$, transfer matrices representing closed-loop responses from disturbance to state and control:
\begin{equation} \label{eq:sf_cl_definitions}
    \xu = \begin{bmatrix} (zI-A-B\K)^{-1} \\ \K(zI-A-B\K)^{-1}\end{bmatrix} \w =: \Phixu \w
\end{equation}

\begin{theorem} \label{thm:sf_main} (Theorem 4 in \cite{Anderson2019}) For system \eqref{eq:plant_x} with state-feedback controller $\u=\K\x$, 
\begin{enumerate}
    \item The affine subspace\footnote{$\frac{1}{z}\RHInf$ is the set of stable, strictly causal LTI transfer matrices.}
\begin{equation} \label{eq:sf_constraint}
    \ZAB\Phixu = I, \hquad \Phix, \Phiu \in \frac{1}{z}\RHInf
\end{equation}
parameterizes all responses $\Phix, \Phiu$ achievable by an internally stabilizing state feedback controller $\K$. This is the \textit{state feedback achievability constraint}.
\item For any $\Phix, \Phiu$ obeying \eqref{eq:sf_constraint}, the controller $\K=\Phiu\Phix^{-1}$, implemented as
\begin{equation}
    \d = \x + (I-z\Phix)\d, \quad \u = z\Phiu\d
\end{equation}
achieves the desired closed-loop response as per \eqref{eq:sf_cl_definitions}. 
\end{enumerate}
\end{theorem}

Any state feedback problem can be cast as \eqref{eq:sls_problem_generic}, where $\Phibf = \Phixu$, and $\Phibf \in \S_a$ is equivalent to $\Phibf$ satisfying \eqref{eq:sf_constraint}.

For the full control problem, we have
\begin{subequations}
\begin{equation}
z\x = A\x + \u + \w
\end{equation}
\begin{equation} \label{eq:plant_y}
\y = C\x + \v
\end{equation}
\end{subequations}
where $\v$ represents measurement noise. Apply linear causal controller $\u = \Lbf\y$ for some transfer matrix $\Lbf$. Define $\Phiw$ and $\Phiv$,  closed-loop responses from disturbance and measurement noise to state:
\begin{equation}
\begin{aligned}
\x &= (zI-A-\Lbf C)^{-1}\w + (zI-A-\Lbf C)^{-1}\Lbf\v \\
   &= \Phiwv \wv
\end{aligned}
\end{equation}
Results from state feedback apply by duality. The \textit{full control achievability constraint} is:
\begin{equation} \label{eq:fc_constraint}
    \Phiwv \ZAC = I, \quad \Phiw, \Phiv \in \frac{1}{z}\RHInf
\end{equation}

An full control problem can be cast as \eqref{eq:sls_problem_generic}, where $\Phibf = \Phiwv$, and $\Phibf \in \S_a$ is equivalent to $\Phibf$ satisfying \eqref{eq:fc_constraint}.

For the output feedback problem defined by \eqref{eq:plant_x} and \eqref{eq:plant_y}, with linear causal controller $\K$ and observer $\Lbf$, we define $\Phixw$, $\Phixv$, $\Phiuw$, $\Phiuv$ -- four closed-loop responses mapping disturbance and measurement noise to state and control:
\begin{equation}
    \xu = \Phiof \wv
\end{equation}

The \textit{output feedback achievability constraint} resembles a combination of the previous achievability constraints \eqref{eq:sf_constraint}, \eqref{eq:fc_constraint}:
\begin{subequations} \label{eq:of_constraint}
\begin{equation} \label{eq:of_col_constraint}
\ZAB \Phiof = \IOHoriz 
\end{equation}
\begin{equation} \label{eq:of_row_constraint}
\Phiof \ZAC = \IOVert
\end{equation}
\begin{equation} \label{eq:of_misc_constraint}
\Phixw, \Phiuw, \Phixv \in \frac{1}{z}\RHInf, \Phiuv \in \RHInf
\end{equation}
\end{subequations}
Detailed derivations are in Section 5 of \cite{Anderson2019}. An output feedback problem can be cast as \eqref{eq:sls_problem_generic}, where $\Phibf = \Phiof$, and $\Phibf \in \S_a$ is equivalent to $\Phibf$ satisfying \eqref{eq:of_constraint}.
\section{Separability and Computation} \label{sec:separability}

We present definitions of separability and compatibility, extending concepts from \cite{Anderson2019}. We describe how separable objectives and constraints in \eqref{eq:sls_problem_generic} translate to distributed computation for both optimal and robust control. 

\textit{Notation:} $\Phibf_{i,j}$ indicates scalar-valued elements of $\Phibf$ that represent influences from the $j^{\text{th}}$ disturbance and/or noise to the $i^{\text{th}}$ node. $\Phibf_{i,:}$ indicates the rows of $\Phibf$ that correspond to the state and/or control at the $i^{\text{th}}$ node. $\Phiset = \{ \Phibf_{i_1, j_1}, \Phibf_{i_2, j_2}, \ldots \}$ indicates some set of elements of $\Phibf$.

\subsection{Separability}
\begin{definition} \textit{(Element-wise separable functional)} Define transfer matrices $\Phibf = \text{argmin} f(\Phibf)$ and $\Psibf$, where $\Psibf_{i,j} = \text{argmin} f_{ij}(\Psibf_{i,j})$. Functional $f(\Phibf)$ is element-wise separable if there exist sub-functionals $f_{ij}$ such that $f(\Phibf) = f(\Psibf)$.
\end{definition}

Relevant element-wise separable functionals include $\htwonorm{\cdot}$ and $\oneinfnorm{\cdot}$ (i.e. maximum absolute element). The latter is the norm associated with $\nu$ robustness \cite{CompanionPaper_Nu}.

\begin{definition} \textit{(Element-wise separable constraint)} Constraint set $\P$ is element-wise separable if there exist sub-constraint sets $\P_{ij}$ such that $\Phibf \in \P \Leftrightarrow \Phibf_{i,j} \in \P_{ij}$.
\end{definition}

Sparsity constraints on $\Phibf$ are element-wise separable.

\begin{definition} \textit{(Row separable functional)}
Define transfer matrices $\Phibf = \text{argmin} f(\Phibf)$ and $\Psibf$, where $\Psibf_{i,:} = \text{argmin} f_i(\Psibf_{i,:})$. Functional $f(\Phibf)$ is row separable if there exist sub-functionals $f_i$ such that $f(\Phibf) = f(\Psibf)$.
\end{definition}

\begin{definition} \textit{(Row separable constraint)} Constraint set $\P$ is row separable if there exist sub-constraint sets $\P_i$ such that $\Phibf \in \P \Leftrightarrow \Phibf_{i,:} \in \P_i$.
\end{definition}

Column separable functionals and constraints are defined analogously. Clearly, any element-wise separable functional or constraint is also both row and column separable. A relevant row separable functional is $\infinfnorm{\cdot}$ (i.e. maximum absolute row sum), which is the norm associated with $\L1$ robustness. A relevant column separable functional is $\oneonenorm{\cdot}$ (i.e. maximum absolute column sum), which is the norm associated with $\LInf$ robustness. It is also trivial to see that constraints of the form $G\Phibf = H$ are column separable, while constraints of the form $\Phibf G = H$ are row separable. Thus, the state feedback achievability constraint is column separable, and the full control achievability constraint is row separable. For the output feedback achievability constraints, \eqref{eq:of_col_constraint} is column separable, \eqref{eq:of_row_constraint} is row separable, and \eqref{eq:of_misc_constraint} is element-wise separable.

\begin{definition} \textit{(Fully separable optimization)} Optimization problem \eqref{eq:sls_problem_generic} is fully row (or column) separable when all objectives and constraints are row (or column) separable.
\end{definition}

\begin{definition} \textit{(Partially separable optimization)} Optimization problem \eqref{eq:sls_problem_generic} is partially separable when all objectives and constraints are either row or column separable, but the overall problem is not fully separable.
\end{definition}

The separability of the SLS state feedback, full control, and output feedback problems are inherently limited by their achievability constraints, since these must always be enforced. The state feedback problem can be fully column separable, but not fully row separable; the opposite is true for the full control problem\footnote{The only exception: if all system matrices are diagonal, then all achievability constraints are element-wise separable.}. The output feedback problem cannot be fully separable, since it includes a mix of column and row separable achievability constraints. We now describe the computational implications of separability.

\subsection{Scalability and Computation}
A fully separable optimization is easily solved via distributed computation. For instance, if problem \eqref{eq:sls_problem_generic} is fully row separable, we can solve the following subproblem:
\begin{equation} \label{eq:sls_subproblem}
    \min_{\Phibf_{i,:}} \quad f_i(\Phibf_{i,:}) \quad
    \textrm{s.t.} \quad \Phibf_{i,:} \in \S_{a_i} \cap \P_i
\end{equation}
at each node $i$ in the system, in parallel.

We often enforce sparsity on $\Phibf$ to constrain inter-node communication to local neighborhoods of size $d$. This translates to $\Phibf_{i,j} = 0 \quad \forall j \notin \N_d(i)$, where $\N_d(i)$ is the set of nodes that are in node $i$'s local neighborhood, and $|\N_d(i)|$ depends on $d$. Then, the size of the decision variable in subproblem \eqref{eq:sls_subproblem} scales with $d$. Each node solves a subproblem in parallel; overall, computational complexity scales with $d$ instead of system size $N$. This is highly beneficial for large systems, where we can choose $d$ much smaller than $N$.

For partially separable problems, we apply the alternating direction method of multipliers (ADMM) \cite{Boyd2010}. We first rewrite \eqref{eq:sls_problem_generic} in terms of row separable and column separable objectives and constraints:
\begin{equation} \label{eq:admm_problem}
    \begin{aligned}
    \min_{\Phibf} \quad & f^{(\row)}(\Phibf) + f^{(\col)}(\Phibf) \\
    \textrm{s.t.} \quad & \Phibf \in \S_a^{(\row)} \cap \S_a^{(\col)} \cap \P^{(\row)} \cap \P^{(\col)} \end{aligned}
\end{equation}
Note that an element-wise separable objective can appear in both $f^{(\row)}$ and $f^{(\col)}$; similarly, an element-wise separable constraint $\P$ can appear in both $\P^{(\row)}$ and $\P^{(\col)}$. We now introduce duplicate variable $\Psibf$ and dual variable $\Lambdabf$. The ADMM algorithm is iterative; for each iteration $k$, we perform the following computations:
\begin{subequations}
\begin{equation} \label{eq:admm_row}
    \begin{aligned}
    \Phibf^{k+1} = \argmin{\Phibf} \quad & f^{(\row)}(\Phibf) + \frac{\gamma}{2}\frobnorm{\Phibf-\Psibf^k+\Lambdabf^k} \\
    \textrm{s.t.} \quad & \Phibf \in \S_a^{(\row)} \cap \P^{(\row)} 
    \end{aligned}
\end{equation}
\begin{equation} \label{eq:admm_col}
    \begin{aligned}
    \Psibf^{k+1} = \argmin{\Psibf} \quad & f^{(\col)}(\Psibf) + \frac{\gamma}{2}\frobnorm{\Phibf^{k+1}-\Psibf+\Lambdabf^k} \\
    \textrm{s.t.} \quad & \Psibf \in \S_a^{(\col)} \cap \P^{(\col)} 
    \end{aligned}
\end{equation}
\begin{equation} \label{eq:admm_lambda}
    \Lambdabf^{k+1} = \Lambdabf^k + \Phibf^{k+1} - \Psibf^{k+1}
\end{equation}
\end{subequations}
ADMM separates \eqref{eq:sls_problem_generic} into a row separable problem \eqref{eq:admm_row} and a column separable problem \eqref{eq:admm_col}, and encourages consensus between $\Phibf$ and $\Psibf$ (i.e. $\Phibf = \Psibf$) via the $\gamma$-weighted objective and \eqref{eq:admm_lambda}. When both $\frobnorm{\Phibf^{k+1}-\Psibf^{k+1}}$ and $\frobnorm{\Phibf^{k+1}-\Phibf^{k}}$ are sufficiently small, the algorithm converges; $\Phibf^k$ is the optimal solution to \eqref{eq:sls_problem_generic}.

Optimizations \eqref{eq:admm_row} and \eqref{eq:admm_col} are fully separable; as described above, they enjoy complexity that scales with local neighborhood size instead of global system size. Additionally, due to sparsity constraints on $\Phibf$ and $\Psibf$, only local communication is required between successive iterations. Thus, partially separable problems also enjoy computational complexity that scale with local neighborhood size instead of global system size \cite{Anderson2019, AmoAlonso2019}. However, partially separable problems require iterations, making them more computationally complex than fully separable problems. The separability of objectives and constraints in \eqref{eq:sls_problem_generic} can also affect convergence rate, as we will describe next.

\begin{definition} \textit{(Compatible sub-functional)} Define $f_\sub$, a functional of $\Phiset$, and $f$, a functional of $\Phibf$. Define transfer matrices $\Phibf$ and $\Psibf$, and let $\Phibf_{i,j} = \Psibf_{i,j} \quad \forall \Phibf_{i,j} \notin \Phiset$. $f_\sub$ is a compatible sub-functional of functional $f$ if minimizing $f_\sub$ is compatible with minimizing $f$, i.e. $f_\sub(\Phiset) \leq f_\sub(\Psiset) \Rightarrow f(\Phibf) \leq f(\Psibf)$
\end{definition}

\begin{definition} \textit{(Compatible sub-constraint)} Constraint set $\P_\sub$ is a compatible sub-constraint of constraint set $\P$ if $\Phibf \in \P \Rightarrow \Phiset \in \P_\sub$ for some set of elements $\Phiset$.
\end{definition}

Compatibility complements separability. If functional $f$ is separable into sub-functionals $f_i$, then each $f_i$ is a compatible sub-functional of $f$; similar arguments apply for constraints.

\begin{definition} \label{defn:balanced} \textit{(Balanced ADMM)}
The partially separable problem \eqref{eq:admm_problem} can be decomposed into row and column sub-functionals and sub-constraints. We say \eqref{eq:admm_problem} is \textit{balanced} if, for any set of elements $\Phiset$ which appear together in a sub-constraint, there exists a matching sub-functional $f_\sub$ that is compatible with $f$, and depends only on $\Phiset$.
\end{definition}

Intuitively, a balanced partially separable problem converges faster than an unbalanced one. For example, consider an output feedback problem with a row separable objective. This is an unbalanced partially separable problem; though all row sub-constraints have a matching sub-objective, none of the column sub-constraints have matching sub-objectives. Thus, we are only able to minimize the objective in the row computation \eqref{eq:admm_row}; this results in slow convergence, and places more burden on consensus between $\Phibf$ and $\Psibf$ than a balanced problem would. More unbalanced examples include a state feedback problem with a row separable objective, or a full control problem with a column separable objective. To balance an output feedback problem, we require an element-wise separable objective $f$.

To summarize: for both fully separable and partially separable problems, computational complexity scales independently of system size. However, partially separable problems require iteration, while fully separable problems do not. For partially separable problems, we prefer a balanced problem to an unbalanced problem due to faster convergence. Thus, element-wise separability (e.g. $\H2$ optimal control, $\nu$ robustness) is desirable for two reasons: firstly, for state feedback and full control, element-wise separable objectives give rise to fully separable problems. Secondly, for output feedback, where ADMM iterations are unavoidable, element-wise separable objectives give rise to balanced problems. Finally, we remark that $\HInf$ robust control problems are not at all separable, and make for highly un-scalable computations; this motivates our use of $\L1$, $\LInf$, and $\nu$ robustness.
\section{Robust Stability} \label{sec:robust_stability}

We present distributed algorithms for structured robust control, using criterion from $\L1$, $\LInf$, and $\nu$ robustness. We consider diagonal nonlinear time-varying (DNLTV) uncertainties $\Delta$ and strictly causal LTI closed-loops $\Phibf$.

\subsection{Robust stability conditions}
\begin{figure}
	\centering
	\begin{minipage}{3cm}
	\includegraphics[width=3cm]{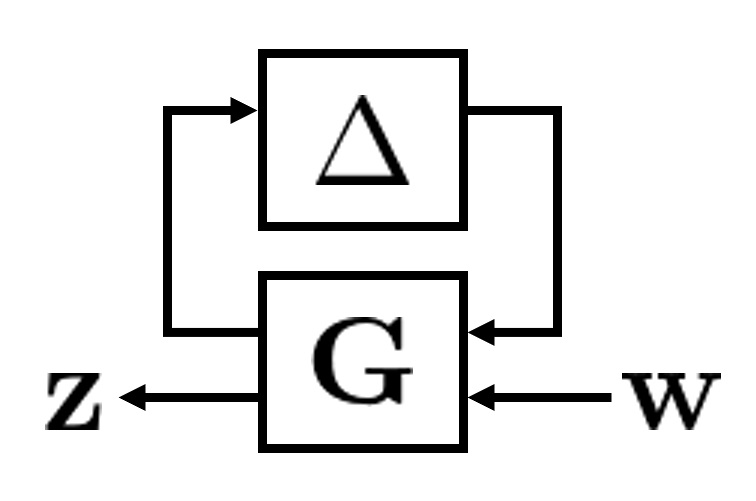}
	\end{minipage}
	\begin{minipage}{4cm}
    \caption{Feedback interconnection of transfer matrix $\G$ and uncertainty $\Delta$. $\G$ is the nominal closed-loop map from disturbance $\w$ to regulated output $\z$.} \label{fig:delta_feedback}
    \end{minipage}
\end{figure}

Let transfer matrix $\G$ map disturbance $\w$ to regulated output $\z$. Generally, $\z$ is a linear function of state $\x$ and input $\u$, i.e. $\G = H\Phibf$ for some constant matrix $H$. Thus, $\G$ is strictly casual and LTI. We construct positive constant magnitude matrix $M = \sum_{p=1}^{\infty} |\G(p)|$, where $\G(p)$ are spectral elements of $\G$, and $|\cdot|$ denotes the element-wise absolute value.
Let $\D$ be the set of positive diagonal matrices:
\begin{equation}
    \D = \{ D \in \mathbb{R}^{n \times n}: D_{ij} = 0 \quad \forall i \neq j, D_{ii} > 0 \quad \forall i\}
\end{equation}

\begin{lemma} \label{lemm:l1_robust} (\textit{$\L1$ robust stability})  The interconnection in Fig. \ref{fig:delta_feedback} is robustly stable for all DNLTV $\Delta$ such that $\infinfnorm{\Delta} < \frac{1}{\beta}$ if and only if $\inf_{D \in \D} \infinfnorm{DMD^{-1}} \leq \beta$. Proof in \cite{Dahleh1993}.
\end{lemma}

\begin{lemma} (\textit{$\LInf$ robust stability}) 
The interconnection in Fig. \ref{fig:delta_feedback} is robustly stable for all DNLTV $\Delta$ such that $\oneonenorm{\Delta} < \frac{1}{\beta}$ if and only if $\inf_{D \in \D} \oneonenorm{DMD^{-1}} \leq \beta$. Proof: equivalent to applying Lemma \ref{lemm:l1_robust} to \ $M^\top$.
\end{lemma}

\begin{lemma} \label{lemm:nu_robust} (\textit{$\nu$ robust stability})
The interconnection in Fig. \ref{fig:delta_feedback} is robustly stable for all DNLTV $\Delta$ such that $\infonenorm{\Delta} < \frac{1}{\beta}$ if $\inf_{D \in \D} \oneinfnorm{DMD^{-1}} \leq \beta$. Additionally, if $\exists D \in \D$ s.t. $DMD^{-1}$ is diagonally maximal, then this condition is both sufficient \textit{and} necessary. Proof: Theorem 4 in \cite{CompanionPaper_Nu}.
\end{lemma}

\begin{definition} \textit{(Diagonally maximal)} Matrix $A \in \mathbb{R}^{n \times n}$ is diagonally maximal if $\exists k$ s.t.  $|A_{kk}|= \max_{i,j}|A_{ij}|$.
\end{definition}

In general, computing the $\infonenorm{\cdot}$ norm is NP-hard; for diagonal $\Delta$, $\infonenorm{\Delta} = \sum_k|\Delta_{kk} |$ \cite{CompanionPaper_Nu}. Let the nominal performance be $\perfnorm{Q\Phibf}$ for some norm $\perfnorm{\cdot}$ and some constant matrix $Q$. Then, the nominal performance and robust stability problem can be posed as follows: \vspace{-1em}

\begin{equation} \label{eq:nonconvex_robstab}
    \begin{aligned}
    \min_{\Phibf, M, D} & \quad \perfnorm{Q\Phibf} + \stabnorm{DMD^{-1}} \\
    \textrm{s.t.} & \quad M = \sum_{p=1}^{T} |H\Phibf(p)|, \hquad
    \Phibf \in \S_a \cap \P, \hquad D \in \D \end{aligned}
\end{equation}
where for ease of computation we assume that $\Phibf$ is finite impulse response (FIR) with horizon $T$. We intentionally leave norms ambiguous; $\stabnorm{\cdot}$ can be $\infinfnorm{\cdot}$, $\oneonenorm{\cdot}$, or $\oneinfnorm{\cdot}$ for $\L1$, $\LInf$, and $\nu$ robust stability, respectively. $\|DMD^{-1}\|$ corresponds to the robust stability margin $\frac{1}{\beta}$; robust stability is guaranteed for all $\Delta$ such that $\|\Delta\| < \frac{1}{\beta}$, for the appropriate norm on $\Delta$. Smaller $\beta$ corresponds to stability guarantees for a larger set of $\Delta$. The nominal performance norm $\perfnorm{\cdot}$ may be different from $\stabnorm{\cdot}$.

\subsection{D-\texorpdfstring{$\Phi$}{} Iteration}
Problem \eqref{eq:nonconvex_robstab} is nonconvex, and does not admit a convex reformulation. Inspired by the D-K iteration method from $\HInf$ robust control \cite{Packard1993}, we adopt an iterative approach. We heuristically minimize \eqref{eq:nonconvex_robstab} by iteratively fixing $D$ and optimizing over $\Phibf$ in the ``$\Phi$ step'', then fixing $\Phibf$ and optimizing (or randomizing) over $D$ in the ``D step''.

\textit{Remark:} Problem \eqref{eq:nonconvex_robstab} poses both nominal performance and robust stability as objectives. If we already know the desired robust stability margin $\betamax^{-1}$, we can omit the stability objective and instead enforce a constraint $\stabnorm{DMD^{-1}} \leq \betamax$, as is done in D-K iteration; Similarly, if we already know the desired nominal performance $\alpha$, we can omit the performance objective and enforce $\perfnorm{Q\Phibf} \leq \alpha$.

The $\Phi$ step solves the following problem:
\begin{equation} \label{eq:phi_step}
    \begin{aligned}
    \Phibf, M = \argmin{\Phibf, M} \quad & \perfnorm{Q\Phibf} \\
    \textrm{s.t.} \quad & \stabnorm{DMD^{-1}} \leq \beta \\
    & M = \sum_{p=1}^{T} |H\Phibf(p)|, \hquad \Phibf \in \S_a \cap \P
    \end{aligned}
\end{equation}
for some fixed value $\beta$ and scaling matrix $D \in \D$. The separability of problem \eqref{eq:phi_step} is characterized by the separability of its objective and constraints. As previously mentioned, $\P$ typically consists of sparsity constraints on $\Phibf$; this constraint is element-wise separable. The separability of other constraints and objectives in \eqref{eq:phi_step} are described below:
\begin{itemize}
    \item If $Q$ is separably diagonal, $\perfnorm{Q\Phibf}$ has the same separability as $\perfnorm{\cdot}$; If not, it is column separable if and only if $\perfnorm{\cdot}$ is column separable. 
    \item $\stabnorm{DMD^{-1}} < \beta$ has the same separability as $\stabnorm{\cdot}$
    \item If $H$ is separably diagonal, $M = \sum_{p=1}^{T} |H\Phibf(p)|$ is element-wise separable; if not, it is column separable. 
    \item $\Phibf \in \S_a$ is column separable for state feedback, row separable for full control, and partially separable for output feedback
\end{itemize}

\begin{definition} \textit{(Separably diagonal)} For state feedback, the product $Q\Phibf$ may be written as $Q_x\Phix + Q_u\Phiu$ for some matrices $Q_x$ and $Q_u$. $Q$ is separably diagonal if both $Q_x$ and $Q_u$ are diagonal matrices. Analogous definitions apply to the full control and output feedback case.
\end{definition}

Table \ref{table:separability_phi_step} summarizes the separability of \eqref{eq:phi_step} for state feedback, full control, and output feedback problems with $\HInf$, $\L1$, $\LInf$, and $\nu$ robustness, where we assume that $Q$ and $H$ are separably diagonal and $\perfnorm{\cdot}$ has the same type of separability as $\stabnorm{\cdot}$. Note that $\HInf$ is not separable for any problem. For state feedback, $\LInf$ and $\nu$ are the preferred stability criteria; for full control, $\L1$ and $\nu$ are the preferred stability criteria. For output feedback, $\nu$ is the only criterion that produces a balanced partially separable problem. Overall, $\nu$ robustness is preferable in all three cases, resulting in either fully separable formulations that require no iterations, or balanced partially separable formulations, which have preferable convergence properties. Though convergence properties vary, the $\Phi$ step \eqref{eq:phi_step} can be computed with complexity that scales with local neighborhood size $d$ instead of global system size for all non-$\HInf$ cases.

\begin{table}
	\caption{Separability of $\Phi$ step of D-$\Phi$ iteration} \vspace{-2em}
	\label{table:separability_phi_step}
	\begin{center}
		\begin{tabular}{|c|c|c|c|}
			\hline
			& State Feedback & Full Control & Output Feedback \\
			\hline
			$\HInf$ & No & No & No \\
			\hline
			$\L1$ & Partial, Unbalanced  & Full & Partial, Unbalanced \\
			\hline
			$\LInf$ & Full & Partial, Unbalanced & Partial, Unbalanced \\
			\hline
			$\nu$ & Full & Full & Partial, Balanced \\
			\hline
		\end{tabular}
	\end{center}
\end{table}

The D step solves the following problem
\begin{equation} \label{eq:d_step_minimize}
    D = \argmin{D} \quad \stabnorm{DMD^{-1}} \quad \textrm{s.t.} \quad D \in \D 
\end{equation}
for some fixed magnitude matrix $M$. For $\nu$ robustness, the minimizing D step \eqref{eq:d_step_minimize} can be recast as a linear program: 
\begin{equation} \label{eq:d_step_minimize_nu}
    \begin{aligned}
    l_i, \eta = \argmin{l_i, \eta} \quad \eta \quad \textrm{s.t.} \quad & \log{M_{i,j}} + l_i - l_j \leq \eta \\ & \forall M_{i,j} \neq 0, \hquad 1 \leq i,j \leq n 
    \end{aligned}
\end{equation}
The optimal solution $D$ can be recovered as $D=\text{diag}(\exp{(l_1)}, \exp{(l_2)}, \ldots \exp{(l_n)})$. Problem \eqref{eq:d_step_minimize_nu} can be distributedly computed using ADMM consensus \cite{Boyd2010}. Let $x_i = \begin{bmatrix} \eta_i \\ L_{j@i} \end{bmatrix}$ be the variable at node $i$, where $\eta_i$ is node $i$'s value of $\eta$, and $L_{j@i}$ is a vector containing $l_{j@i}$: node $i$'s values of $l_j$ for all $j \in \N_d(i)$. The goal is for all nodes to reach consensus on $\eta$ and $l_i, 1 \leq i \leq N$. We introduce dual variable $y_i$ and averaging variable $\bar{x}_i = \begin{bmatrix} \bar{\eta}_i \\ \bar{L}_{j@i} \end{bmatrix}$. For each iteration $k$, node $i$ performs the following computations:
\begin{subequations} \label{eq:admm_consensus_nu}
\begin{equation} \label{eq:admm_consensus_1}
    \begin{aligned}
    x_i^{k+1} & = \argmin{x_i} \quad \eta_i + (y_i^k)^\top(x_i - \bar{x}_i^k) + \gamma\|x_i - \bar{x}_i^k \|_2^2 \\
    & \textrm{s.t.} \quad \forall j \in \N_d(i), M_{i,j} + l_{i@i} - l_{j@i} \leq \eta_i
    \end{aligned}
\end{equation}
\begin{equation} \label{eq:admm_consensus_2}
    \begin{aligned}
    \bar{\eta}_i^{k+1} & = \frac{1}{|\N_d(i)|}\sum_{j \in \N_d(i)} \eta_j \\
    \bar{l}_i & = \frac{1}{|\N_d(i)|}\sum_{j \in \N_d(i)} l_{i@j}, \quad \bar{l}_{j@i}^{k+1} = \bar{l}_j, \forall j \in \N_d(i)
    \end{aligned}
\end{equation}
\begin{equation} \label{eq:admm_consensus_3}
    y_i^{k+1} = y_i^k + \frac{\gamma}{2}(x_i^{k+1} - \bar{x}_i^{k+1})
\end{equation}
\end{subequations}
where $\gamma$ is a user-determined parameter, and iterations stop when consensus is reached, i.e. differences between relevant variables are sufficiently small. The size of optimization variable $x_i$ depends only on local neighborhood size $d$; thus, the complexity of \eqref{eq:admm_consensus_1} scales independently of global system size. Computation \eqref{eq:admm_consensus_2} requires communication, but only within the local neighborhood. Also, by the definition of $\Phibf$ and the construction of $M$, $M_{i,j} = 0 \quad \forall j \notin \N_d(i)$. Thus, for a fully connected system, solving \eqref{eq:admm_consensus_nu} is equivalent to solving \eqref{eq:d_step_minimize_nu}. Additionally, consensus problems are balanced as per Definition \ref{defn:balanced}, so \eqref{eq:admm_consensus_nu} converges relatively quickly.

For $\L1$ robustness, \eqref{eq:d_step_minimize} is solved by setting $D = \text{diag}(v_1, v_2, \ldots v_n)^{-1}$, where $v$ is the eigenvector corresponding to the largest-magnitude eigenvalue of $M$ \cite{Dahleh1993}. This computation does not lend itself to scalable distributed computation. To ameliorate this, we propose an alternative formulation that randomizes instead of minimizing over $D$: 
\begin{equation} \label{eq:d_step_randomize}
    D = \argmin{D} \quad 0 \quad \textrm{s.t.} \quad \stabnorm{DMD^{-1}} \leq \beta
\end{equation}
The randomizing formulation lends itself to distributed computation. Also, we remark that \eqref{eq:d_step_minimize} can be solved by iteratively solving \eqref{eq:d_step_randomize} to search for the lowest feasible $\beta$. 

Define vectors $\vec{D} = \begin{bmatrix} D_{11} & D_{22} & \ldots & D_{nn} \end{bmatrix}^\top$ and $\vec{D}^{-1} = \begin{bmatrix} D_{11}^{-1} & D_{22}^{-1} & \ldots & D_{nn}^{-1} \end{bmatrix}^\top$. Then, we can rewrite constraint $\stabnorm{DMD^{-1}} \leq \beta$ as $M\vec{D}^{-1} \leq \beta \vec{D}^{-1}$ for $\L1$ stability, and $M^\top \vec{D} \leq \beta \vec{D}$ for $\LInf$ stability. Then, problem \eqref{eq:d_step_randomize} can be formulated as a scalable distributed ADMM consensus problem using similar techniques as \eqref{eq:d_step_minimize_nu}.

Both versions of the D step (\eqref{eq:d_step_minimize} and \eqref{eq:d_step_randomize}) for D-$\Phi$ iteration are simpler than the D step in D-K iteration \cite{Packard1993}, which requires a somewhat involved frequency fitting process. Also, all separable versions of the proposed D step are less computationally intensive than the $\Phi$ step \eqref{eq:phi_step}, since the decision variable in the D step is much smaller. Table \ref{table:scalability_D_step} summarizes the scalability of different versions of the D step for different robustness criteria. `Minimize' refers to solving \eqref{eq:d_step_minimize} directly; `Iteratively Minimize' refers to solving \eqref{eq:d_step_minimize} by iteratively solving \eqref{eq:d_step_randomize} to search for the lowest feasible $\beta$; `Randomize' refers to solving \eqref{eq:d_step_minimize}. \cmark \hquad indicates that we can use scalable distributed computation, and \xmark \hquad indicates that no scalable distributed formulation is available; by scalable, we mean complexity that scales independently of global system size. For iterative minimization, there is the obvious caveat of iterations incurring additional computational time; however, for $\L1$ and $\LInf$, iterative minimization is more scalable than direct minimization. Additionally, we show in the next section that algorithms using the randomizing D step perform similarly as algorithms using the minimizing D step; thus, iterative minimization may be unnecessary. Overall, $\nu$ robustness appears to be preferable for scalability purposes for both the $\Phi$ step and D step.

\begin{table}
	\caption{Scalability of D step of D-$\Phi$ iteration} \vspace{-1em}
	\label{table:scalability_D_step}
	\begin{center}
		\begin{tabular}{|c|c|c|c|}
		    \hline
		    & Minimize & Iteratively Minimize & Randomize \\
			\hline
			$\HInf$ & \xmark & \xmark & \xmark \\
			\hline
			$\L1$ & \xmark & \cmark & \cmark \\
			\hline
			$\LInf$ & \xmark & \cmark & \cmark \\
            \hline
			$\nu$ & \cmark & \cmark & \cmark \\
            \hline            
		\end{tabular}
	\end{center}
\end{table}

We now present two algorithms for D-$\Phi$ iteration. Algorithm \ref{alg:d_phi_iters_minimizing} is based on minimizing over $D$ \eqref{eq:d_step_minimize}, while Algorithm \ref{alg:d_phi_iters_randomizing} is based on randomizing over $D$ \eqref{eq:d_step_randomize}. Both algorithms compute the controller $\Phibf$ which achieves optimal nominal performance for some desired robust stability margin $\betamax^{-1}$. \vspace{-0.5em}

\begin{algorithm}
\caption{D-$\Phi$ iteration with minimizing D step} \label{alg:d_phi_iters_minimizing}
\begin{algorithmic}[1]
\Statex \textbf{input:} $\betastep > 0$, $\betamax > 0$ \quad \textbf{output:} $\Phibf$
\State Initialize $\beta^{k=0} \leftarrow \infty$, $k \leftarrow 1$
\State Set $\beta^k \leftarrow \beta^{k-1} - \betastep$. Solve \eqref{eq:phi_step} to obtain $\Phibf^k$, $M^k$
\Statex \textbf{if} \eqref{eq:phi_step} is infeasible\textbf{:}
\quad \textbf{return} $\Phibf^{k-1}$, $\beta^{k-1}$
\State Solve \eqref{eq:d_step_minimize} to obtain $D$. Set $\beta^k \leftarrow \stabnorm{DM^kD^{-1}}$
\Statex \textbf{if} $\beta^k \leq \betamax$\textbf{:} \quad \textbf{return $\Phibf^k$, $\beta^k$}
\State Set $k \leftarrow k+1$ and return to step 2
\end{algorithmic}
\end{algorithm} \vspace{-0.5em}

\begin{algorithm}
\caption{D-$\Phi$ iteration with randomizing D step} \label{alg:d_phi_iters_randomizing}
\begin{algorithmic}[1]
\Statex \textbf{input:} $\betastep > 0$, $\betamax > 0$ \quad \textbf{output:} $\Phibf$
\State Initialize $\beta^{k=0} \leftarrow \infty$, $k \leftarrow 1$, $D \leftarrow I$
\State Set $\beta^k \leftarrow \beta^{k-1} - \betastep$. Solve \eqref{eq:phi_step} to obtain $\Phibf^k$, $M^k$
\Statex \textbf{if} $k=1$\textbf{:} \quad Set $\beta^k \leftarrow \stabnorm{DM^kD^{-1}}$
\Statex \textbf{if} \eqref{eq:phi_step} is infeasible\textbf{:} \quad Solve \eqref{eq:d_step_randomize}
\Statex \quad \quad \textbf{if} \eqref{eq:d_step_randomize} is infeasible\textbf{:} \quad \textbf{return} $\Phibf^{k-1}$, $\beta^{k-1}$
\Statex \quad \quad \textbf{else:} \quad Solve \eqref{eq:phi_step} to obtain $\Phibf^k$, $M^k$
\Statex \textbf{if} $\beta^k \leq \betamax$\textbf{:} \quad \textbf{return $\Phibf^k$, $\beta^k$}
\State Solve \eqref{eq:d_step_randomize} to obtain $D$
\State Set $k \leftarrow k+1$ and return to step 2
\end{algorithmic}
\end{algorithm}

In Algorithm \ref{alg:d_phi_iters_minimizing}, we alternate between minimizing over $\Phibf$ and minimizing over $D$, and stop when no more progress can be made or when $\betamax$ is attained. No initial guess of $D$ is needed; at iteration $k=1$, $\beta^k = \infty$, and the $\stabnorm{DMD^{-1}} \leq \beta$ constraint in \eqref{eq:phi_step} of step 2 is trivially satisfied. In Algorithm \ref{alg:d_phi_iters_randomizing},  we alternate between minimizing $\Phibf$ and \textit{randomizing} $D$. There are two main departures from Algorithm \ref{alg:d_phi_iters_minimizing} due to the use of the randomizing D step:
\begin{enumerate}
    \item An initial guess of $D$ is required to generate $\beta^{k=1}$, which is then used as an input to the randomizing D step. $D=I$ is a natural choice, although we may also randomize or minimize over the initial $D$.
    \item In step 2, when we cannot find a new $\Phibf$ to make progress on $\beta$, instead of stopping, attempt to find a new $D$ to make progress on $\beta$. If such a $D$ can be found, find the new optimal $\Phibf$, then continue iterating.
\end{enumerate}

Parameter $\betastep$ appears in both algorithms, and indicates the minimal robust stability margin improvement per step. For both algorithms, computational complexity is dominated by the $\Phi$ step problem \eqref{eq:phi_step} and D step problem \eqref{eq:d_step_minimize} or \eqref{eq:d_step_randomize}. All of these problems can be distributedly computed, and all enjoy complexity that scales independently of global system size; thus, the complexity of the overall algorithm also scales independently of global system size.
\section{Simulations} \label{sec:simulations}
We use a ring of size $N=10$ with spectral radius 3:
\begin{equation}
\begin{aligned}
    x_i(t+1) & = A_{i,i-1}x_{i-1}(t) + A_{i,i+1}x_{i+1}(t) + u_i(t) \\
    & \text{for} \quad i = 2 \ldots N-1, \\
    x_1(t+1) & = A_{1,N}x_{N}(t) + A_{1,2}x_{2}(t) + u_1(t), \\
    x_N(t+1) & = A_{N,N-1}x_{N-1}(t) + A_{N,1}x_{1}(t) + u_N(t)
\end{aligned}
\end{equation}
Nonzero elements of system matrix $A$ are randomly drawn. We focus on the state feedback case, with an LQR nominal performance objective with state penalty $Q_x = I$ and control penalty $Q_u = 50I$. The regulated output is $z_i = x_i + u_i$. 

\begin{figure}
	\centering
	\includegraphics[width=\columnwidth]{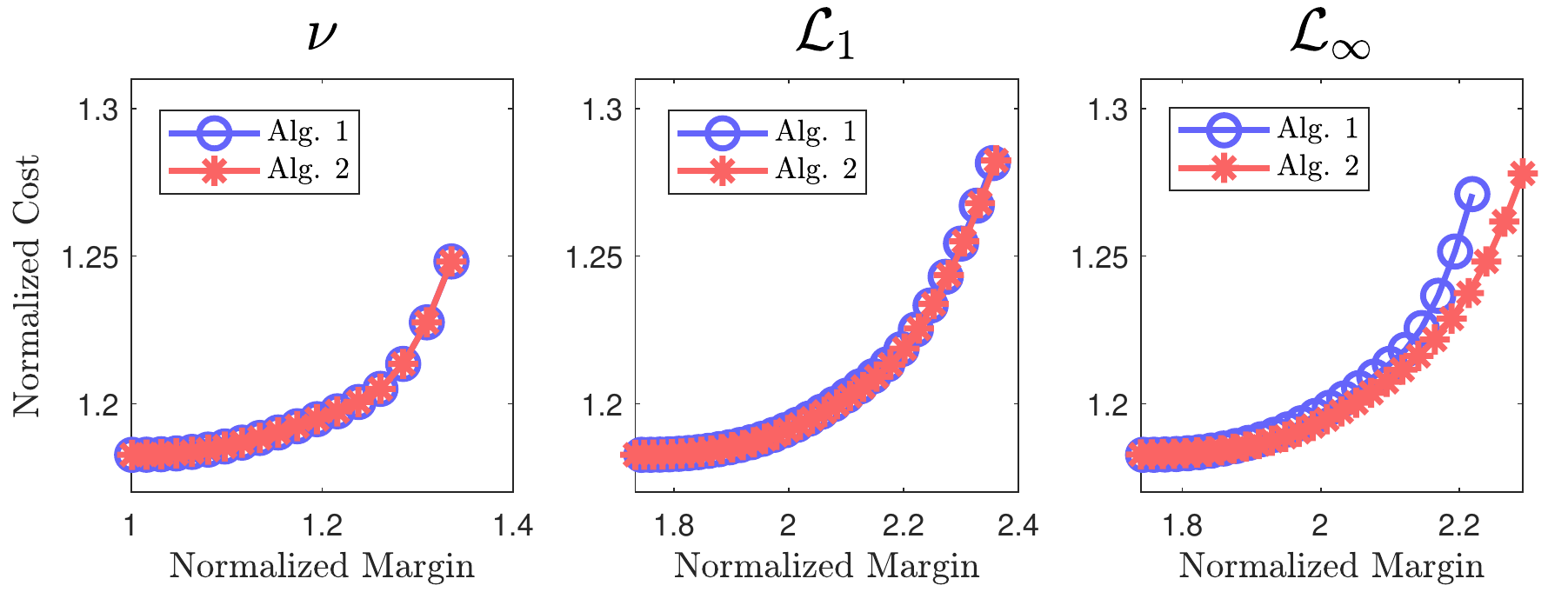}
	\vspace{-1.5em}
	\caption{D-$\Phi$ iteration results for $\nu$, $\L1$, and $\LInf$ robust stability. Algorithm \ref{alg:d_phi_iters_minimizing} and Algorithm \ref{alg:d_phi_iters_randomizing} perform similarly. The controller with maximum stability margin for $\nu$ is found in 18 iterations for both algorithms; for $\L1$, in 30 iterations for Algorithm \ref{alg:d_phi_iters_minimizing} and 35 iterations for Algorithm \ref{alg:d_phi_iters_randomizing}; for $\LInf$, for $\L1$, in 25 iterations for Algorithm \ref{alg:d_phi_iters_minimizing} and 35 iterations for Algorithm \ref{alg:d_phi_iters_randomizing}.}
	\label{fig:random_chain}
\end{figure}

We constrain $\Phibf$ to be FIR with horizon size $T=30$. We also constrain $\Phibf$ to be sparse, such that each node is only allowed to communicate to neighbors and neighbors of neighbors. We run Algorithm \ref{alg:d_phi_iters_minimizing} and \ref{alg:d_phi_iters_randomizing} with $\betastep = 0.05$ and varying values of $\betamax$, and compare the results to the optimal LQR solution, which is unconstrained by local communication. Results are shown in Fig. \ref{fig:random_chain}; cost and robust stability margin are normalized against the LQR solution.

Both algorithms start with the standard SLS solution; in each plot, the point on the bottom left corresponds to this solution. The performance suboptimality of the SLS solution arises purely from communication constraints. There appears to be a tradeoff between nominal cost and stability margin; as the margin increases, so does the cost. Note that as the number of iterations increase, robust stability margin $\frac{1}{\beta}$ can only increase; thus, the number of iterations required increases as the desired margin $\betamax^{-1}$ increases. Algorithm \ref{alg:d_phi_iters_minimizing} and Algorithm \ref{alg:d_phi_iters_randomizing} perform similarly, though Algorithm \ref{alg:d_phi_iters_randomizing} appears to require more iterations. Overall, these simulations show that the two proposed algorithms for D-$\Phi$ iteration are viable and give good margins and nominal performance, as desired. Simulation results may be reproduced using the relevant scripts in the SLS-MATLAB toolbox \cite{Li2019_SLSMatlab}.
\section{Future work} \label{sec:future_work}
This paper provides a novel algorithm for distributed structured robust control. Historically, structured robust control theory, particularly $\HInf$, was not formulated with distributed synthesis in mind. However, distributed and scalable synthesis is becoming increasingly relevant for modern large-scale systems. Immediate avenues of future work include \textit{1)} incorporating block diagonal uncertainty, and \textit{2)} incorporating robust \textit{performance} in addition to robust stability.

\bibliography{refs}

\begin{thebibliography}{10}
\providecommand{\url}[1]{#1}
\csname url@samestyle\endcsname
\providecommand{\newblock}{\relax}
\providecommand{\bibinfo}[2]{#2}
\providecommand{\BIBentrySTDinterwordspacing}{\spaceskip=0pt\relax}
\providecommand{\BIBentryALTinterwordstretchfactor}{4}
\providecommand{\BIBentryALTinterwordspacing}{\spaceskip=\fontdimen2\font plus
\BIBentryALTinterwordstretchfactor\fontdimen3\font minus
  \fontdimen4\font\relax}
\providecommand{\BIBforeignlanguage}[2]{{%
\expandafter\ifx\csname l@#1\endcsname\relax
\typeout{** WARNING: IEEEtran.bst: No hyphenation pattern has been}%
\typeout{** loaded for the language `#1'. Using the pattern for}%
\typeout{** the default language instead.}%
\else
\language=\csname l@#1\endcsname
\fi
#2}}
\providecommand{\BIBdecl}{\relax}
\BIBdecl

\bibitem{Dahleh1993}
M.~A. Dahleh and M.~H. Khammash, ``Controller design for plants with structured
  uncertainty,'' \emph{Automatica}, vol.~29, pp. 37--56, 1993.

\bibitem{Packard1993}
A.~Packard and J.~C. Doyle, ``The complex structured singular value,''
  \emph{Automatica}, vol.~29, no.~1, pp. 71--109, 1993.

\bibitem{Zhou1998}
K.~Zhou and J.~C. Doyle, \emph{Essentials of robust control}.\hskip 1em plus
  0.5em minus 0.4em\relax Prentice Hall New Jersey, 1998.

\bibitem{Bamieh2002}
B.~Bamieh, F.~Paganini, and M.~A. Dahleh, ``{Distributed Control of Spatially
  Invariant Systems},'' \emph{IEEE Transactions on Automatic Control}, vol.~47,
  no.~7, pp. 1091--1107, 2002.

\bibitem{Lidstrom2017}
C.~Lidstrom, R.~Pates, and A.~Rantzer, ``{H-infinity optimal distributed
  control in discrete time},'' in \emph{IEEE CDC}, 2017, pp. 3525--3530.

\bibitem{Matni2020}
N.~Matni and A.~A. Sarma, ``Robust performance guarantees for system level
  synthesis,'' in \emph{IEEE ACC}, 2020, pp. 779--786.

\bibitem{Anderson2019}
J.~Anderson, J.~C. Doyle, S.~H. Low, and N.~Matni, ``{System level
  synthesis},'' \emph{Annual Reviews in Control}, vol.~47, pp. 364--393, 2019.

\bibitem{CompanionPaper_Nu}
O.~Kjellqvist and J.~C. Doyle, ``$\nu$-analysis: A new notion of robustness for
  large systems with structured uncertainties,'' to appear in IEEE CDC 2022.

\bibitem{Boyd2010}
S.~Boyd, N.~Parikh, E.~Chu, B.~Peleato, and J.~Eckstein, ``Distributed
  optimization and statistical learning via the alternating direction method of
  multipliers,'' \emph{Foundations and Trends in Machine Learning}, vol.~3,
  no.~1, pp. 1--122, 2011.

\bibitem{AmoAlonso2019}
C.~{Amo Alonso} and N.~Matni, ``{Distributed and Localized Closed Loop Model
  Predictive Control via System Level Synthesis},'' in \emph{IEEE CDC}, 2020,
  pp. 5598--5605.

\bibitem{Li2019_SLSMatlab}
\BIBentryALTinterwordspacing
J.~S. Li, ``{SLS-MATLAB}: Matlab toolbox for system level synthesis,'' 2019.
  [Online]. Available: \url{https://github.com/sls-caltech/sls-code}
\BIBentrySTDinterwordspacing

\end{thebibliography}
\bibliographystyle{IEEEtran}

\end{document}